\documentstyle[leqno]{bcp}
\def\LaTeX{L\kern -.36em\raise .3ex\hbox{\sc a}\kern -.15em T\kern
-.1667em%
\lower .7ex\hbox{E}\kern -.125em X}

\mathchardef\zd="710E  
\mathchardef\zG="7000  
\mathchardef\zL="7003  
\mathchardef\zl="7115  

\newcommand{\be}{\begin{equation}}
\newcommand{\ee}{\end{equation}}
\newcommand{\ra}{\rightarrow}

\newcommand{\bea}{\begin{eqnarray}}
\newcommand{\eea}{\end{eqnarray}}
\newcommand{\beas}{\begin{eqnarray*}}
\newcommand{\eeas}{\end{eqnarray*}}
\newcommand{\1}{{\bf 1}}
\newcommand{\k}{{\bf k}}
\newcommand{\we}{\wedge}
\newcommand{\nn}{\nonumber}

\newcommand{\pa}{\partial}
\newcommand{\ti}{\times}
\newcommand{\A}{{\cal A}}

\begin{document}
\mathclass{Primary 17B65 81R10; Secondary 17B63 37J99 53D99 70H99.}
\thanks{Research of the second author supported by PRIN SINTESI.}
\abbrevauthors{J. Grabowski and G. Marmo}
\abbrevtitle{Binary operations}

\title{Binary operations \\ in classical and quantum mechanics}

\author{Janusz\ Grabowski}
\address{Mathematical Institute, Polish Academy of Sciences\\
ul. \'Sniadeckich 8, P.O. Box 137, 00-950 Warszawa, Poland\\
E-mail: jagrab@mimuw.edu.pl}

\author{Giuseppe\ Marmo}
\address{Dipartimento di Scienze Fisiche, Universit\`a Federico II di
Napoli\\ and\\ INFN, Sezione di Napoli\\ Complesso Universitario di Monte
Sant'Angelo\\ Via Cintia, 80126 Napoli, Italy\\
E-mail:marmo@na.infn.it}

\maketitlebcp

\abstract{Binary operations on algebras of observables are studied
in the quantum as well as in the classical case. It is shown that certain
natural compatibility conditions with the associative product imply
properties which are usually additionally required.}

\section*{1. Introduction.}
It was an observation of P.~A.~M. Dirac, when considering the quantum
Poisson bracket of observables in foundations of Quantum Mechanics, that
the Leibniz rule
\beas
[A,B_1B_2]&=&B_1[A,B_2]+[A,B_1]B_2, \cr
[B_1B_2,A]&=&B_1[B_2,A]+[B_1,A]B_2,
\eeas
is sufficient to determine the bracket. The Leibniz rule simply tells
us that fixing an argument in the bracket we get a derivation, i.e. an
infinitesimal automorphism, of the algebra of observables.
This shows that, in fact, most of the axioms of the Lie
bracket for the corresponding operator algebra are superfluous
and one can insist only on the proper behaviour with respect to the
composition (Leibniz rule).
\par
In the classical case, however, the Leibniz rule is not enough to
determine proper classical brackets, since every operation on the algebra
$C^\infty(M)$ of smooth functions on a manifold $M$, associated with the
contraction with any contravariant
2-tensor, will do. We have to impose a version of the Jacobi identity:
$$
[x,[y,z]]=[[x,y],z]+[y,[x,z]].
$$
However, as it was observed in [GM], the Jacobi identity and the Leibniz
rule for a bracket on $C^\infty(M)$ imply the skew-symmetry, i.e. we are
dealing with a Poisson bracket.
\par
Note that just the Jacobi identity (without the skew-symmetry assumption)
is the axiom of a generalization of a Lie algebra proposed by J.-L.~Loday
[Lo]. Such non-skew-symmetric Lie algebras do exist and give rise to a
well-defined (co)homology theory. The result of [GM] just shows that such
Loday algebras are not possible (except for the skew-symmetric ones, of
course) on algebras of functions when the Leibniz rule is required.
The natural question arises, if we can break the skew-symmetry by passing
to more general structures like Jacobi brackets or local Lie algebras in
the sense of A.~A.~Kirillov [Ki].
\par
Our main result shows that this is impossible. More precisely, we show
that, assuming the differentiability of the bracket and the Jacobi
identity, we get the skew-symmetry automaticaly, so that no local Loday
brackets do exist except for the well-known skew-symmetric ones.
The differentiability condition may be wiewed as a compatibility condition
with the associative product in the algebra.
These observations seem interesting to us in the context of finding a
minimal set of requirements for the structures of classical and
quantum mechanics.

\section*{2. Commutator brackets in associative rings.}

It was already observed by Dirac ([Di], p.86) that the Leibniz
rule defines the commutator bracket in the algebra of quantum observables
uniquely up to a constant factor.
In fact, the following holds for any associative ring (cf. [Di]).

\th{Proposition}{1.}{If $[\cdot,\cdot]_0$ is a bilinear operation in an
associative ring $\A$ for which the Leibniz rule is
satisfied (with respect to both arguments):
\bea\label{lr}
[A,B_1B_2]_0&=&B_1[A,B_2]_0+[A,B_1]_0B_2, \cr
[B_1B_2,A]_0&=&B_1[B_2,A]_0+[B_1,A]_0B_2,
\eea
then
\be\label{1}
[A_1,B_1][A_2,B_2]_0=[A_1,B_1]_0[A_2,B_2],
\ee
for all $A_i,B_i\in\A$, where $[A,B]=AB-BA$ is the commutator bracket.}

\Proof Using the Leibniz rule first for the second argument and then for
the first one, we get
\beas
[A_1A_2,B_1B_2]_0&=&B_1[A_1A_2,B_2]_0+[A_1A_2,B_1]_0B_2=\cr
&&B_1A_1[A_2,B_2]_0+B_1[A_1,B_2]_0A_2+A_1[A_2,B_1]_0B_2+[A_1,B_1]_0A_2B_2.
\eeas
Doing the same calculations in the reversed order we get
$$
[A_1A_2,B_1B_2]_0=
A_1B_1[A_2,B_2]_0+B_1[A_1,B_2]_0A_2+A_1[A_2,B_1]_0B_2+[A_1,B_1]_0B_2A_2.
$$
Hence
$$A_1B_1[A_2,B_2]_0+[A_1,B_1]_0B_2A_2=B_1A_1[A_2,B_2]_0+[A_1,B_1]_0A_2B_2
$$
and the proposition follows
\sq
\th{Theorem}{1.}{If $\A$ is a unital associative ring which is {\it
strongly non-commutative},
that is the two-sided associative ideal of $\A$ generated by the derived
algebra $\A'=\mathop{span}\{[A,B]:A,B\in\A\}$ equals $\A$, then, under the
assumptions of Proposition 1, we get
\be
[A,B]_0=C[A,B]
\ee
for some central element $C\in\A$.
\par
In particular, every bilinear operator $[\cdot,\cdot]_0$ on $\A$ which
satisfies the Leibniz rule is a Lie bracket, i.e. it is skew-symmetric and
satisfies the Jacobi identity. If, moreover, $\A$ is an algebra over a
commutative ring $\k$ and the center of $\A$ is
trivial, i.e. equals $\k\1$, then $[A,B]_0=\zl[A,B]$ for certain
$\zl\in\k$.}

\Proof By assumption, the unit $\1\in\A$ can be written in the form
\be
\1=\sum_iC_i[A_i,B_i]+\sum_j[A_j',B_j']C_j'.
\ee
From (\ref{1}) we get then
\be
\left(\sum_iC_i[A_i,B_i]\right)[A,B]_0=
\left(\sum_iC_i[A_i,B_i]_0\right)[A,B].
\ee
Since both brackets satisfy the Leibniz rule, we have
\bea\label{2}
&&[A_j',B_j']C_j'[A,B]_0=[A_j',B_j'C_j'][A,B]_0-B_j[A_j',C_j'][A,B]_0=\cr
&&[A_j',B_j'C_j']_0[A,B]-B_j[A_j',C_j']_0[A,B]=
[A_j',B_j']_0C_j'[A,B],
\eea
so that
\beas
[A,B]_0&=&\left(\sum_iC_i[A_i,B_i]+\sum_j[A_j',B_j']C_j'\right)[A,B]_0=\cr
&&\left(\sum_iC_i[A_i,B_i]_0+\sum_j[A_j',B_j']_0C_j'\right)[A,B]=C[A,B],
\eeas
where $C$ is a fixed element of $\A$ not depending on $A,B\in\A$.
\par
To show that $C$ is central, rewrite (\ref{2}) in the form
\be\nn
[A',B']XC[A,B]=C[A',B']X[A,B],
\ee
so that
\be\label{3}
[[A',B']X,C]\cdot[A,B]Y=0,
\ee
where $A,B,A',B',X,Y$ are arbitrary elements of $\A$. Since every
one-sided associative ideal of $\A$ containing the derived Lie algebra
$\A'$ is, due to the Leibniz rule, two-sided and thus equals $\A$, we get
from (\ref{3}) $[\A,C]\A=0$, whence $[\A,C]=0$ and the theorem follows.
\sq

\remar{Remark\ {1.}\ }{Obviously, if $\A$ contains elements $P,Q$
satisfying the canonical commutation rules $PQ-QP=\1$, then the
assumption of the above theorem is satisfied automaticaly. This is
exactly the Dirac's case.}

\th{Corollary}{1.}{Let $\A$ be the algebra $\mathop{gl}(n,\k)$ of $n\ti
n$-matrices with coefficients in a commuatative unital ring $\k$. Then
every binary operation $[\cdot,\cdot]_0$ on $\A$, satisfying the Leibniz
rule (\ref{lr}), is of the form
\be\label{cc}
[A,B]_0=\zl(AB-BA).
\ee
for certain $\zl\in\k$.}

\Proof It is easy to see that the derived Lie algebra $\A'$ is the Lie
algebra $\mathop{sl}(n,\k)$ of trace-less matrices. Of course,
$\mathop{gl}(n,\k)\ne\mathop{sl}(n,\k)$,
but $\mathop{sl}(n,\k)$ contains invertible matrices, so that the
associative ideal generated by $\A'$ is the whole $\A$. Due to theorem 1,
$[A,B]_0=C(AB-BA)$ for some central element $C\in\mathop{gl}(n,\k)$.
Since, as easily seen, the center of $\mathop{gl}(n,\k)$ is trivial, we
get (\ref{cc}).
\sq

\section*{3. Differential Loday brackets for commutative algebras.}

The important bra\-ckets of classical mechanics, like the Poisson
bracket on a symplectic manifold or the Lagrange bracket on a contact
manifold, were generalized by A.~A.~Kirillov [Ki] to {\it local
Lie algebra} brackets on one-dimensional vector bundles over a manifold
$M$, that is to Lie
brackets given by local operators. The fundamental fact discovered in
[Ki] is that these operators have to be of the first order and then,
locally, they reduce to the conformally symplectic Poisson and Lagrange
brackets on the leaves of the corresponding generalized foliation of $M$.
For the trivial bundle, i.e. for the algebra $C^\infty(M)$ of functions on
$M$, the local brackets reduce to the so called {\it Jacobi brackets}
associated with the corresponding {\it Jacobi structures} on the manifold
$M$ (cf.[Li]).
\par
A purely algebraic version of the Kirillov's result has been proved in
[Gr], theorems 4.2 and 4.4, where Lie brackets on associative
commutative algebras with no nilpotents, given by bidifferential operators,
have been considered.

\smallskip\noindent
On the other hand, J.-L.~Loday (cf. [Lo1]), while studying relations
between Hochschild and cyclic homology in the search for
obstructions to the periodicity of algebraic K-theory, discovered that
one can skip the skew-symmetry assumption in the definition of a Lie algebra,
still having a possibility to define appropriate (co)homology
(see [Lo1,LP] and [Lo], Chapter 10.6). His Jacobi
identity for such structures was formally the same as the classical Jacobi
identity in the form
\be\label{JI}
[x,[y,z]]=[[x,y],z]+[y,[x,z]].
\ee
This time, however, this is no longer equivalent to
\be\label{JI1}
[[x,y],z]=[[x,z],y]+[x,[y,z]],
\ee
nor to
\be
[x,[y,z]]+[y,[z,x]]+[z,[x,y]]=0,
\ee
since we have no skew-symmetry. Loday called such structures {\it
Leibniz algebras} but, since we have already associated the name of
Leibniz with the Leibniz rule, we shall call them {\it Loday algebras}.
This is in accordance with the terminology of [KS], where
analogous structures in the graded case are defined. Of course, there is
no particular reason not to define Loday algebras by means of (\ref{JI1})
instead of (\ref{JI}) (and in
fact, it was the original Loday definition), but both categories are
equivalent via transposition of arguments. Similarly, for associative
algebras we can obtain associated algebras by transposing arguments, but
in this case we still get associative algebras.
\par
The natural question arises, whether local Loday algebras, different from
the local Lie algebras, do exist. In [GM] it has been proved
that no new $n$-ary Loday-Poisson algebras and, in principle, no new Loday
algebroids (in the sense of [GU]) are possible.
Here we will prove that local Loday algebras on sections of
one-dimensional vector bundles must be skew-symmetric, that is, they
reduce
to local Lie algebras of Kirillov.
We will start with a general algebraic result on this subject.
\par
Let now $\A$ be an associative commutative algebra with the unit $\1$ over
a field $\k$ of characteristic 0.
Recall that a (linear) {\it differential operator of order} $\le n$ on
$\A$  is a $\k$-linear operator $D:\A\ra\A$ such that $\zd(x)^{n+1}D=0$
for all $x\in\A$, where $\zd(x)$ is the commutator with the multiplication
by $x$:
\be
(\zd(x)D)(y)=D(xy)-xD(y).
\ee
This is the same as to say that $\zd(x_1)\cdots\zd(x_{n+1})D=0$ for any
$x_1,\dots,x_{n+1}\in\A$. Note that a zero-order differential operator $D$
is just the multiplication by $D(\1)$ and that $\zd(x)$ acts as a
derivation of the associative algebra $\mathop{End}_\k(\A)$ of linear
operators on $\A$:
\be
\zd(x)(D_1D_2)=(\zd(x)D_1)D_2+D_1(\zd(x)D_2).
\ee
For multilinear operators we define anologously
\be
\zd_{y_i}(x)D(y_1,\dots,y_p)=D(y_1,\dots,xy_i,\dots,y_p)-xD(y_1,\dots,y_p)
\ee
the corresponding derivations with respect to the $i$'th variable and call
the multilinear operator $D$ {\it being of order} $\le n$ if
$\zd_{y_i}(x)^{n+1}D=0$ for all $x\in\A$ and all $i$, i.e. $D$ is of order
$\le n$ with respect to each variable separately. 
This means that fixing $(p-1)$ arguments we get a differential operator of
order $\le n$. One can also consider multilinear operators such that,
fixing $(p-1)$ arguments, we get a differential operator with some order
depending (and possibly unbounded) on what we have fixed. In this paper,
however, by a multilinear differential operator we mean an operator of
order $n$ for some $n$. Note also that the differentials $\zd_{y_i}(x)$ and 
$\zd_{y_j}(u)$ commute.

\defin{Definition}{1. {\it A differential Loday bracket} on $\A$ is a
Loday algebra bracket on $\A$, i.e. a bracket satisfying (\ref{JI}), given
by a bidifferential operator of certain order $n$.}

\th{Proposition}{2.}{If $\A$ has no nontrivial nilpotent elements, then
every differential Loday bracket on $\A$ is of order $\le 1$.}

\Proof Let $D:\A\ti\A\ni(x,y)\ra[x,y]\in\A$ be a differential Loday
bracket. Denote $D_x=[x,\cdot]$, $D_x^r=[\cdot,x]$, and let $n$ be the
order of $D$ with respect to the second argument, i.e. the maximum of
orders of $D_x$ for all $x\in\A$. Then, $\zd(w)^{n}D_x$ is the
multiplication by $\zd(w)^nD_x(\1)$  for all
$w,x\in\A$ and $\zd(w)^nD_x(\1)$ is different from zero for some $w,x\in\A$.
Let now $k$ be the maximum of orders of $\zd(w)^nD_x(\1)$ with respect to $x$.
Obviously, $k$ is not greater that the order $m$ of $D$ with respect to
the first variable. In other words,
\be\label{0}
\zd_{x}(u)^p(\zd(w)^qD_x)=\zd_0(u)^p\zd_1(w)^qD(x,\cdot)=0,
\ee
(where $\zd_0,\zd_1$ are the corresponding differentials with respect to
the first and the second variable, respectively), if $q>n$, or $q=n$ and
$p>k$, and
\be
\zd_{x}(u)^k(\zd(w)^nD_x)=\zd_0(u)^k\zd_1(w)^nD(x,\cdot)
\ee
is the multiplication by
\be
\zd_0(u)^k\zd_1(w)^nD(\1,\1)
\ee
which is different from zero for certain $u=u_0$ and $w=w_0$.
Note that since we do not assume that
the bracket is skew-symmetric, we do not have $D_x+D_x^r=0$ and that $n$
is the order of $D$ with respect to the first argument.
We claim first that $n\le 1$.
\par
Suppose the contrary. Then $2n-1>n$ and
\bea
&&\zd_y(v)^k\zd_x(u)^{k+1}\zd_z(w)^{2n-1}(D(D(x,y),z)+D(y,D(x,z))-
D(x,D(y,z))=\cr
&&(k+1){2n-1\choose n}\left(\zd_0(u)^k\zd_1(w)^nD(\1,\1)
\cdot\zd_0^k(v)\zd_1(u)\zd_1(w)^{n-1}D(\1,\1)\right),\label{4}
\eea
where we identify the right-hand-side term with the corresponding
zero-order differential operator.
Indeed, $\zd(w)^{2n-1}$ vanishes on the first summand $D_{[x,y]}$.
For the rest we get
\bea
&&\zd(w)^{2n-1}[D_y,D_x]_c(z)=\cr
&&{2n-1\choose n}\left([\zd(w)^{n}D_y,\zd(w)^{n-1}D_x]_c +
[\zd(w)^{n-1}D_y,\zd(w)^{n}D_x]_c\right)(z),
\eea
where the brackets are the commutator brackets.
Since $\zd(w)^{n}D_y$ is a multiplication by $\zd(w)^{n}D_y(\1)$, the
first summand is just $\zd(\zd(w)^{n}D_y(\1))\zd(w)^{n-1}D_x(z)$, so it
vanishes after applying $\zd_x(u)^{k+1}$ in view of (\ref{0}). From the
second summand, after applying $\zd_x(u)^{k+1}$, we get clearly
\bea
&&\zd_x(u)^{k+1}(\zd(w)^{n-1}D_y\circ\zd(w)^{n}D_x)(z)=\cr
&&(k+1)\left(\zd_0(u)^k\zd_1(w)^nD(\1,\1)\cdot\zd_1(u)\zd_1(w)^{n-1}
D_y(z)\right).
\eea
Applying now $\zd_y(v)^k$ we get (\ref{4}). Since
\be\label{j}
D(D(x,y),z)+D(y,D(x,z))-D(x,D(y,z))=0
\ee
due to the Jacobi identity, we get
\be\label{5}
\zd_0(u)^k\zd_1(w)^nD(\1,\1)\cdot\zd_0^k(v)\zd_1(u)\zd_1(w)^{n-1}
D(\1,\1)=0.
\ee
Putting now $u:=u+tw$ in (\ref{5}) with $t$ varying through
$\k$ and passing to the first derivative with respect to $t$ (the
coefficient by $t$) we get
\bea\label{6}
0&=&k\zd_0(u)^{k-1}\zd_0(w)\zd_1(w)^nD(\1,\1)\cdot\zd_0^k(v)\zd_1(u)
\zd_1(w)^{n-1}D(\1,\1)+\cr
&&\zd_0(u)^k\zd_1(w)^nD(\1,\1)\cdot\zd_0^k(v)\zd_1(w)^{n}D(\1,\1).
\eea
Multiplicating both sides of (\ref{6}) by $\zd_0(u)^k\zd_1(w)^nD(\1,\1)$
and using (\ref{5}) we get, after putting $u=v$,
\be
(\zd_0(u)^k\zd_1(w)^nD(\1,\1))^3=0,
\ee
thus
\be
\zd_0(u)^k\zd_1(w)^nD(\1,\1)=0
\ee
since there are no nontrivial nilpotents in $\A$. But the latter is
different from zero for $u=u_0$, $w=w_0$; a contradiction.
\par
Thus $n\le 1$ and we will finish with showing that $m$ -- the order of $D$
with respect to the first argument -- is also $\le 1$.
As above, suppose the contrary and let $s$ be the maximal order
of $\zd(w)^mD_y^r(\1)$ with respect to $y$, so that there are
$w_0,u_0\in\A$ such that
\be\label{c}
\zd_0(w_0)^m\zd_1(u_0)^sD(\1,\1)\ne 0.
\ee
We get from (\ref{j})
\bea\label{n}
0&=&\zd_z(u)^s\zd_y(u)^s\zd_x(w)^{2m}\left(D(D(x,y),z)+D(y,D(x,z))-
D(x,D(y,z))\right)=\cr
&&{2m\choose m}\left(\zd_1(u)^s\zd_0(w)^mD(\1,\1)\right)^2
\eea
which contradicts (\ref{c}). Indeed, since $2m>m>1$ and
\be
D(D(x,y),z)+D(y,D(x,z))-D(x,D(y,z))=(D_z^rD_y^r+D_yD_z^r-D_{[y,z]})(x),
\ee
$\zd_x(w)^{2m}$ vanishes on $D_{[y,z]}$ ($2m>m$) and on $D_yD_z^r$ (we
already know that $D_y$ is of order $\le 1$), so for the middle term of
(\ref{n}) we get
\bea
&&\zd_z(u)^s\zd_y(u)^s\zd_x(w)^{2m}\left(D_z^r\cdot D_y^r\right)(x)=\cr
&&x{2m\choose m}\zd_z(u)^s\zd_y(u)^s
\left(\zd(w)^m(D_z^r)\cdot\zd(w)^m(D_y^r)\right)=\cr
&&xyz{2m\choose m}\left(\zd_1(u)^s\zd_0(w)^mD(\1,\1)\right)^2.
\eea
\sq

\remar{Remark\ {2.}\ }{The bi-orders $(k,n)$ and $(m,s)$ correspond to
bidifferential parts of $D$ of maximal orders with respect to
the anti-lexicographical and the lexicographical orderings, respectively.
We could say that they are the corresponding bi-symbols of this
bidifferential operator. Note, however, that it is not true in general
algebraic case that differential operators of higher orders are
polynomials in derivatives as for differential operators on $C^\infty(M)$
(cf. [Gr], remark 3.2).
This is only for first-order differential operators on $\A$ that we have
the decomposition $\mathop{Der}(\A)\oplus\A$ into the direct sum of
derivations and zero-order operators (multiplications by elements of
$\A$).}

\remar{Remark\ {4.}\ }{The assumption concerning the nilpotent elements is
essential. For instance, if $\A$ is freely generated by elements $x,y$
with the constraint $x^2=0$, then
\be
D(u,v)=xu\frac{\pa^n}{\pa y^n}(v)
\ee
is a differential Loday bracket of order $n$.}

\th{Theorem}{2.}{If $\A$ has no nontrivial nilpotent elements, then
every differential Loday bracket on $\A$ is a standard
skew-symmetric Jacobi bracket given by
\be
[x,y]=\zL(x,y)+x\zG(y)-y\zG(x),
\ee
where $\zL$ and $\zG$ are, respectively, a biderivation  and a derivation
on $\A$ which satisfy the compatibility conditions
\bea
[\zG,\zL]_{NR}&=&0;\cr
[\zL,\zL]_{NR}&=&-2\zL\we\zG,
\eea
for $[\cdot,\cdot]_{NR}$ being the Nijenhuis-Richardson bracket of
skew-multilinear maps.}

\remar{Remark\ {4.}\ } {The above theorem generalizes the theorem 4.4 in
[Gr]. The Nijenhuis-Richardson graded Lie algebra bracket [NR] (see also
[Gr]) is in this algebraic case the analog of the Schouten-Nijenhuis
bracket of multivector fields.}

\Proof In view of proposition 2, we can assume that the bracket
$D=[\cdot,\cdot]$ is of first-order and, according to [Gr], theorem 4.4.,
it is sufficient to prove that the bracket is skew-symmetric.
The bracket being of the first order satisfies the generalized Leibniz
rule:
\bea\label{e}
[x,yz]&=&y[x,z]+[x,y]z-yz[x,1];\cr
[yz,x]&=&y[z,x]+[y,x]z-yz[1,x].
\eea
Indeed, $D_x$ is of the first order, so $\zd(z)D_x$ is the multiplication
by $\zd(z)D_x(\1)$. But $\zd(z)D_x(y)=[x,yz]-z[x,y]$, so
\be
[x,yz]-z[x,y]=(\zd(z)D_x(\1))y=([x,z]-z[x,1])y.
\ee
One proves analogously the second equation.
The equations (\ref{e}) mean that the operators $R_x$ and $L_y$ defined by
\bea
L_x(y)&=&[x,y]-[x,\1]y,\cr
R_y(x)&=&[x,y]-[\1,y]x,
\eea
are derivations of $\A$. They are just obtained from the decompositions of
the corresponding first-order differential operators mentioned in remark
4.
Moreover,
\bea
(\zd(y)D_x)(z)&=&L_x(y)z,\cr
(\zd(y)D_x^r)(z)&=&R_x(y)z.
\eea
From the Jacobi identity (\ref{JI}) we get immediately that
\be\label{l}
[[x,x],z]=0
\ee
and, using the Jacobi identity once more, that
\be\label{r}
[[y,[x,x]],z]=0
\ee
for all $x,y,z\in\A$. The linearized version of (\ref{l}) reads
\be\label{li}
[[x,y]+[y,x],z]=0.
\ee
From (\ref{li}) we get easily
\bea
0&=&\zd_x(u)^2([[x,y]+[y,x],z])=x\zd(u)^2(D_z^r(D_y+D_y^r))(1)=\cr
&&2xR_z(u)(L_y(u)+R_y(u)),
\eea
so
\be\label{imp1}
R_z(u)(L_y(u)+R_y(u))=0.
\ee
In view of (\ref{l}), $D_{[x,x]}=0$, so that $L_{[x,x]}=0$. But putting
$y=z=[x,x]$ in (\ref{imp1}) we get $(R_{[x,x]}(u))^2=0$, so that
$R_{[x,x]}=0$ and hence
\be\label{imp2}
[z,[x,x]]=z[\1,[x,x]].
\ee
Putting now $y=[1,[x,x]]$ in (\ref{imp2}) we get, in view of (\ref{r}),
\be
([1,[x,x]])^2=[[1,[x,x]],[x,x]]=0,
\ee
so that $[1,[x,x]]=0$ and finally
\be\label{new}
[z,[x,x]]=0.
\ee
Now, using (\ref{new}) instead of (\ref{l}) we get, similarly to
(\ref{imp1}),
\be
L_z(u)(L_y(u)+R_y(u))=0
\ee
which, together with (\ref{imp1}), yields
\be
(L_z(u)+R_z(u))(L_y(u)+R_y(u))=0.
\ee
The latter, after putting $z=y$, implies cleary $L_y+R_y=0$.
But $(L_y+R_y)(x)=0$ means
\be\label{w}
[x,y]+[y,x]=x([y,\1]+[\1,y]).
\ee
Since the right-hand side of (\ref{w}) is symmetric with respect to $x$
and $y$,
\be
x([y,\1]+[\1,y])=y([x,\1]+[\1,x])
\ee
which implies
\be\label{p}
[y,\1]+[\1,y]=2y[\1,\1].
\ee
Putting in (\ref{p}) $y=[\1,\1]$ we get, due to (\ref{l}),
$[\1,\1]^3=0$, so that
\be
[y,\1]+[\1,y]=0
\ee
and the skew-symmetry $[x,y]+[y,x]=0$ follows from (\ref{w}).
\sq

\th{Corollary}{2.}{Every local binary operation on sections of a
one-dimensional bundle over a manifold $M$ which satisfies the Jacobi
identity (\ref{JI}) is skew-symmetric, i.e. it is a local Lie algebra
bracket.\par
In particular, every such operation on the algebra $C^\infty(M)$ is a
Jacobi bracket associated with a Jacobi structure on $M$.}

\Proof Since the operation is local, we can reduce locally to the algebra
of smooth functions and, due to locality of the operation, to a
bidifferential operator. In the algebra of smooth functions there are no
nontrivial nilpotents, so, according to theorem 2, the operation is of
first order and skew-symmetric.
\sq


\references{GM}
{\item{[Di]} {\sc P. A. M. Dirac}, {\it The Principles of Quantum
Mechanics}\/, Oxford University Press, Oxford 1958.

\item{[Gr]}
{\sc J. Grabowski},
{\it Abstract Jacobi and Poisson structures}\/,
J.~Geom.~Phys.  9 (1992), 45--73.

\item{[GM]}
{\sc J. Grabowski}  and {\sc G. Marmo},
{\it Non-antisymmetric version of Nambu-Poisson and algebroid brackets}\/,
J.~Phys~A: Math.Gen. 34 (2001), 3803--3809.

\item{[GU]}{\sc J. Grabowski} and {\sc P. Urba\'nski},
{\it Algebroids -- general differential calculi on vector bundles}\/,
J.~Geom.~Phys.  31 (1999), 111-141.

\item{[Ki]}{\sc A. A. Kirillov},
{\it Local Lie algebras}\/, Russ.~Math.~Surv. 31 (1976), 55--75.

\item{[KS]} {\sc Y. Kosmann-Schwarzbach},
{\it From Poisson algebras to Gerstenhaber algebras}\/,
\\ Ann.~Inst.~Fourier 46 (1996), 1243--1274.

\item{[Li]} {\sc A. Lichnerowicz},
{\it Les vari\'et\'es de Jacobi et leurs alg\`ebres de Lie associ\'ees}\/,
J.~Math. Pures Appl. 57 (1978), 453--488.

\item{[Lo]}
{\sc J.-L. Loday},
{\it Une version non commutative des alg\`ebres de Lie: les alg\`ebres
de Leibniz}\/, Ann.~Inst.~Fourier 37 (1993), 269--93.

\item{[Lo1]} {\sc J.-L. Loday},
{\it Cyclic Homology}\/,
Springer Verlag, Berlin 1992.

\item{[LP]}
{\sc J.-L. Loday  and T. Pirashvili},
{\it Universal enveloping algebras of Leibniz algebras and
(co)homology}\/, Math.~Annalen 296 (1993), 569--572.

\item{[NR]} {\sc A. Nijenhuis} and {\sc R. Richardson},
{\it Deformation of Lie algebra brackets}\/,\\
J.~Math.~Mech. 17 (1967), 89--105.
}
\end{document}